\documentclass[12pt]{amsart}
\usepackage{amsmath, amsthm, amssymb}
\usepackage{enumerate}

\theoremstyle{plain}

\theoremstyle{plain}
\newtheorem{thm}{Theorem}
\newtheorem{cor}[thm]{Corollary}
\newtheorem{lem}[thm]{Lemma}
\newtheorem{prop}[thm]{Proposition}

\theoremstyle{definition}

\newtheorem{defi}{Definition}
\newtheorem{exa}{Example}
\newtheorem{exas}{Examples}
\newtheorem{const}{Construction}
\newtheorem{sit}{}

\theoremstyle{remark}

\newtheorem{rem}{Remark}
\newtheorem{rems}{Remarks}

\usepackage{enumerate}
\usepackage{amssymb}
\input{diagrams.tex}

\newcommand{\brem}{\begin{rem}}
\newcommand{\brems}{\begin{rems}}
\newcommand{\erem}{\end{rem}}
\newcommand{\erems}{\end{rems}}
\newcommand{\bexa}{\begin{exa}}
\newcommand{\bexas}{\begin{exas}}
\newcommand{\eexa}{\end{exa}}
\newcommand{\eexas}{\end{exas}}
\newcommand{\bdefi}{\begin{defi}}
\newcommand{\edefi}{\end{defi}}
\newcommand{\bdefis}{\begin{defis}}
\newcommand{\edefis}{\end{defis}}
\newcommand{\bcor}{\begin{cor}}
\newcommand{\ecor}{\end{cor}}
\newcommand{\blem}{\begin{lem}}
\newcommand{\elem}{\end{lem}}
\newcommand{\bconv}{\begin{conv}}
\newcommand{\econv}{\end{conv}}
\newcommand{\bconj}{\begin{conj}}
\newcommand{\econj}{\end{conj}}
\newcommand{\bprop}{\begin{prop}}
\newcommand{\eprop}{\end{prop}}
\newcommand{\bthm}{\begin{thm}}
\newcommand{\ethm}{\end{thm}}
\newcommand{\bnota}{\begin{nota}}
\newcommand{\enota}{\end{nota}}
\newcommand{\bconst}{\begin{const}}
\newcommand{\econst}{\end{const}}
\newcommand{\bsit}{\begin{sit}}
\newcommand{\esit}{\end{sit}}
\newcommand{\be}{\begin{eqnarray}}
\newcommand{\ee}{\end{eqnarray}}
\newcommand{\bproof}{\begin{proof}}
\newcommand{\eproof}{\end{proof}}
\def\ba{\begin{array}}
\def\ea{\end{array}}
\newcommand{\no}{\noindent}

\def\cO{{\mathcal O}}

\newcommand{\A}{{\mathbb A}}
\newcommand{\PP}{{\mathbb P}}
\newcommand{\R}{{\mathbb R}}
\newcommand{\C}{{\mathbb C}}
\newcommand{\Q}{{\mathbb Q}}
\newcommand{\Z}{{\mathbb Z}}
\newcommand{\N}{{\mathbb N}}

\newcommand{\G}{{\Gamma}}

\begin{document}



\title[Affine lines on $\Q$-homology planes and group actions]
{Affine lines on $\Q$-homology planes and group actions}

\author{Mikhail Zaidenberg}
\thanks{
\mbox{\hspace{11pt}}{\it 1991 Mathematics Subject Classification}:
14R05, 14R20, 14J50.\\
\mbox{\hspace{11pt}}{\it Key words}: $\Q$-acyclic surface,
$\Q$-homology plane}
\email{zaidenbe@ujf-grenoble.fr}
\address{Universit\'e
Grenoble I\\ Institut Fourier\\ UMR 5582 CNRS-UJF\\ BP 74\\ 38402
St.\ Martin d'H\`eres c\'edex\\ France}






\begin{abstract} This note is a supplement to the papers \cite{KiKo}
and \cite{GMMR}.
We show the role of group actions in
classification of affine lines on $\Q$-homology planes.
\end{abstract}

\maketitle


\section*{Introduction}

This note\footnote{To appear in Transformation Groups 11:4, 2006;
available on line at:\\
https://webmail.math.cnrs.fr/cgi-bin/nph-revues.cgi/010110A/http
/www.springerlink.com/content/2823v3t286257h72/fulltext.pdf} is a
supplement to the papers \cite{KiKo} and \cite{GMMR}. Our aim is
to shed a light on the role of group actions in classification of
affine lines on $\Q$-homology planes with logarithmic Kodaira
dimension $-\infty$. This enables us to strengthen certain results
in {\it loc. sit.} (see Section 1).

Let us fix terminology. It is usual \cite[Ch. 3, \S 4]{Mi} to call
a smooth $\Q$-acyclic ($\Z$-acyclic, respectively) surface over
$\C$ a {\it $\Q$-homology plane} (a {\it homology plane},
respectively). By Fujita's Lemma \cite[2.5]{Fu} such a surface is
necessarily affine. Likewise we call a {\it homology line} an
irreducible affine curve $\Gamma$ with Euler characteristic
$e(\Gamma)=1$. So $\G$ is homeomorphic to $\R^2$ and its
normalization is isomorphic to $\A^1=\A^1_\C$. A smooth curve
isomorphic to $\A^1$ will be called an {\it affine line}.
Following \cite{Mi} we let $\A^1_*=\A^1\setminus\{0\}$. As usual
$\bar k$ stands for logarithmic Kodaira dimension.

{\bf Acknowledgements:} {\footnotesize This research was done
during a visit of the author to the Max Planck Institute of
Mathematics in Bonn. He thanks this institution for a generous
support and excellent working conditions. Our thanks also to
Shulim Kaliman and to a referee of the 'Transformation groups" for
useful editorial comments.}

\section{Main results}

\bthm\label{theo} Let $X$ be a $\Q$-homology plane and $\Gamma$
a homology line on $X$.
Then the following hold.
\begin{enumerate}
[(a)] \item Suppose that $\bar k(X\setminus\Gamma)=-\infty$. Then
$\G$ is either an orbit of an effective $\C_+$-action on $X$ or a
connected component of the fixed point set of such an action.
Anyhow $\G\simeq\A^1$ is a fiber component of the corresponding
orbit map (an $\A^1$-ruling) $\pi:X\to\A^1$. \item Suppose that
$\bar k(X\setminus\Gamma)\ge 0$. Suppose further that
$\Gamma\simeq\A^1$ and $\bar k(X)=-\infty$. Then $\Gamma$ is an
orbit closure of an effective hyperbolic $\C^*$-action on $X$.
Moreover $X$ admits an effective action of a semidirect product
$G=\C^* \ltimes\C_+$ with an open orbit $U$. The orbit map
$X\to\A^1$ of the induced $\C_+$-action defines an $\A^1$-ruling
on $X$ with a unique multiple fiber say $\G'\simeq\A^1$ such that
$\G$ and $\G'$ meet at one point transversally and $U=X\setminus
\G'\simeq \A^1\times\A^1_*$. Furthermore this $\C_+$-action moves
$\G$. Consequently there exists a continuous family of affine
lines $\G_t$ on $X$ with the same properties as $\G$. \item
Suppose that  $\Gamma$ is singular. Then $X\simeq\A^2$ and $\bar k
(X\setminus \G)=1$. Moreover\footnote{This is due to the
Lin-Zaidenberg Theorem \cite[Ch. 3, \S 3]{LiZa, Mi}.} there is an
isomorphism $X\simeq\A^2$ sending $\Gamma$ to a curve $V(x^k-y^l)$
with coprime $k,l\ge 2$. Consequently $\G$ is an orbit closure of
an elliptic $\C^*$-action on $X$.
\end{enumerate}\ethm

\no We indicate below a proof of the theorem. The cases
(a), (b) and (c) are proven in
Sections 2, 3 and 4, respectively. Besides,
in cases (a) and (b) we
provide in Lemmas \ref{constr1} and \ref{dpd},
respectively,
a description of the pairs $(X,\G)$
satisfying
their assumptions.
The assertion of (b) follows from Theorem 1.1
in \cite{KiKo},
cf. also Theorem 3.10 in \cite{GMMR}.
In the case of a $\Z$-homology plane
(c) was established in \cite{Za}; the proof for
a $\Q$-homology plane is similar.
This gives a strengthening of Theorem 1.3 in
\cite{KiKo}.

The cases (a)-(c) of Theorem \ref{theo}
do not exhaust all the possibilities for the
pair $(X,\G)$ as above. To complete the picture let us summarize
some known facts,
see e.g. \cite[Ch. 3, \S 4]{Za, GuPa, Mi}
and the references therein.

\bthm\label{corMT} We let as before $X$ be a $\Q$-homology plane
and $\G\subseteq X$ a homology line.
If $\G$ is singular then $(X,\G)$ is as in Theorem \ref{theo}(c).
Suppose further that $\G$ is
smooth i.e. is an affine line. Then
$\bar k(X)\le\bar k(X\setminus\G)\le 1$\footnote{See
\cite[Ch.2, Theorem 6.7.1]{Mi}.}
and one of the following cases occurs.
\begin{enumerate}
[(a)] \item $\left(\bar k(X),\, \bar
k(X\setminus\G)\right)=(-\infty,-\infty)$ and $(X,\G)$ is as in
Theorem \ref{theo}(a) that is, $\G$ is of fiber type and
$X\setminus\G$ carries a family of disjoint affine
lines\footnote{See also Lemma \ref{constr1} below.}. \item
$\left(\bar k(X),\, \bar k(X\setminus\G)\right) =(-\infty,0)$ or
$(-\infty,1)$ and $(X,\G)$ is as in Theorem
\ref{theo}(b)\footnote{The both possibilities actually occur, see
the Construction in Section 3 and also Lemma \ref{dpd}.}. \item
$\left(\bar k(X),\, \bar k(X\setminus\G)\right)=(0,0)$ and either
$X$ is not NC minimal or $X$ is one of the Fujita's surfaces
$H[-k,k]$ ($k\ge 1$)\footnote{We refer e.g. to \cite[Ch. 3,
4.4.1-4.4.2]{Fu, GuPa, Mi} for definitions.}. Anyhow $\G$ is a
unique affine line on $X$ unless $X=H[-1,1]$. \item $\left(\bar
k(X),\, \bar k(X\setminus\G)\right)=(0,1)$, $X=H[-1,1]$ and there
are exactly two affine lines, say, $\G_0$ and $\G_1=\G$ on $X$.
These lines meet transversally in two distinct points, moreover
$\bar k(X\setminus\G_0)=0$ and $\bar k(X\setminus\G_1)=1$. \item
$\left(\bar k(X),\, \bar k(X\setminus\G)\right)=(1,1)$, there is a
unique $\A^1_*$-fibration on $X$ and $\G$ is a fiber component of
its degenerate fiber\footnote{The same conclusions hold also in
case (c) if $X$ is not NC-minimal \cite{GuPa}.}. There can be at
most one further affine line on $X$, which is then another
component of this same degenerate fiber, and these components meet
transversally in one point.
\end{enumerate}
\ethm

\brem\label{newrem} Let $X$ be a $\Z$-homology plane.
By \cite{Fu} then $\bar k (X)\neq 0$. By \cite{Za} (supplement)
$\bar k (X)=1$ if and only if there exists a unique homology
(in fact, affine) line on $X$.
 \erem

\renewcommand\thesubsection{\arabic{subsection}}
\section{$\Q$-homology planes with an $\A^1$-ruling}
These occur to be smooth affine surfaces with
$\A^1$-rulings $X\to\A^1$ which possess only irreducible
degenerate fibers. They were studied in details e.g. in \cite[4.14]{Fu},
\cite{Be}, \cite{Fi}, \cite[\S 4]{FlZa1}. See also \cite[Ch. 3, 4.3.1]{Mi} for
a brief summary\footnote{We note \cite{Be} that $\pi_1(X)$ is a free product
of cyclic groups, namely, $\pi_1(X)\cong *_j \Z/m_j\Z$, where $(m_j)_j$ is
the sequence of multiplicities of degenerate fibers, and so
$H_1(X;\Z)\cong \bigoplus_j \Z/m_j\Z$.}.
In Lemma \ref{constr1} below we show that every
$\A^1$-ruling $\pi:X\to \A^1$ on a $\Q$-homology plane $X$
can be obtained starting from a standard linear $\A^1$-ruling $\A^2\to\A^1$
and replacing several fibers by multiple fibers via a procedure called
in \cite{FlZa1} a {\it comb attachment}.
More precisely, this replacement goes as follows.

{\it  Attaching combs}. On the quadric $\PP^1\times \PP^1$ with a
$\PP^1$-ruling $\pi_0={\rm pr}_1: \PP^1\times \PP^1\to \PP^1$ we
fix a finite set of points $\{A_j\}$, $j=1,\ldots,n$ ($n\ge 0$) in
different fibers $F_j=\{t_j\}\times \PP^1$ of $\pi_0$. We fix
further a sequence $\sigma: V\to \PP^1\times \PP^1$ of blowups
with centers at the points $A_j$ and infinitesimally near points.
Letting $\bar\pi: V\to \PP^1$ be the induced $\PP^1$-ruling, we
suppose that $\bar\pi$ enjoys the following properties:
\begin{enumerate}
[$\bullet$] \item the center
of every blowup over $A_j$ except for the first one belongs
to the exceptional
($-1$)-curve of the previous blowup;
\item $D_\infty\cdot E_j=0$ $\forall j=1,\ldots,n$, where $D_{\infty}$
is the proper transform
in $V$ of the section $\PP^1\times\{\infty\}$ of ${\rm pr}_1$ and
$E_j$ is
the last ($-1$)-curve in the fiber $\bar\pi^{-1}(t_j)$.
\item $E_j$ is a tip of the dual graph
of the fiber $\bar\pi^{-1}(t_j)$.
\end{enumerate}
Under these assumptions the dual graph as above is a comb, with
all vertices of degree $\le 3$. Let
$F_{\infty}=\bar\pi^{-1}(t_\infty)\subset V$ be a fiber over an
extra point $t_{\infty}\in \PP^1\setminus\{t_1,\ldots,t_n\}$ and
$E\subseteq V$ be the reduced exceptional divisor of $\sigma: V\to
\PP^1\times \PP^1$. We consider the open surface $X=V\setminus D$,
where $D=F_{\infty}+D_{\infty}+E+\sum_{j=1}^n (F_j'-E_j)$ and
$F_j'$ is the proper transform  of $F_j$ in $V$. Then $\bar\pi : V
\to \PP^1$ restricts to an $\A^1$-ruling $\pi : X \to \A^1$ with
only irreducible fibers; all fibers of $\pi$ are reduced except
possibly the fibers $\pi^{-1}(t_j)=E_j\cap X$.

The following lemma is well known, see e.g.
\cite[Proposition 4.9]{FlZa1}.

\blem\label{constr1} Under the notation as above the surface $X$
is a $\Q$-homology plane. Moreover, every $\Q$-homology plane $X$
with an $\A^1$-ruling $\pi : X \to \A^1$ arises in this way. \elem

\bproof Let $X$ be constructed as above.
By the Suzuki formula \cite{Suz, Za, Gu}, $e(X)=1$
and so the equality $b_2=b_1+b_3$ for the Betti numbers of $X$ holds.
Thus $X$ is $\Q$-acyclic
if and only if $b_2=0$
or equivalently,
if ${\rm Pic} (D)\otimes \Q$ generates
${\rm Pic} (V)\otimes \Q$, see \cite[Ch. 3,
4.2.1]{Mi}. The latter
is easily seen to be the case in our construction. The
first assertion follows now
by Fujita's Lemma
\cite[2.5]{Fu}.

As for the second one, given an $\A^1$-ruling $\pi : X \to
\A^1$ on a $\Q$-homology plane $X$ it extends
to a pseudominimal
$\PP^1$-ruling $\pi:V\to \PP^1$ on
a  smooth
completion $V$ of $X$ with an SNC boundary divisor $D$.
The pseudominimality
means that none of the
$(-1)$-curves in $D-D_\infty$, where $D_\infty$ is the horizontal
component of $D$, can be
contracted without loosing the SNC property, see \cite[3.4]{Za}.
Since $e(X)=1$ all fibers of $\pi : X \to
\A^1$ are irreducible. We let $\bar\pi^{-1}(t_j)$,
$j=1,\ldots,n$ ($n\ge 0$)
be the degenerate fibers of $\bar\pi$ and $E_j$
be the component of the fiber
$\bar\pi^{-1}(t_j)$
such that $E_j\cap
X=\pi^{-1}(t_j)\simeq \A^1$.
By the pseudominimality assumption, $E_j$
is the only $(-1)$-curve
in the fiber
$\bar\pi^{-1}(t_j)$.
Therefore $V$ is obtained from a Hirzebruch surface
$\Sigma_m$
by blowing up process which enjoys
the properties of a comb attachment.
Performing, if necessary,
elementary transformations
in the fiber $F_\infty$ we
may assume that $\bar D_\infty^2=0$, where $\bar D_\infty$ is the image
of $D_\infty$ in $\Sigma_m$ and so,
$\Sigma_m=\Sigma_0=\PP^1\times\PP^1$.
\eproof

\brem\label{afmod} Every surface $X$ as considered above can
actually be obtained from affine plane $\A^2$ via a suitable
affine modification that is \cite[\S 1]{KaZa}, by blowing up with
center in a zero dimensional subscheme $V(I)$ of $\A^2$ located on
a principal divisor $D$ and deleting the proper transform of $D$.
Indeed $X$ contains a cylinder $U\times \A^1$, where
$U\subseteq\A^1$ is a Zariski open subset, see \cite{MiSu} or
\cite[Ch. 3, 1.3.2]{Mi}. The canonical projections of $U\times
\A^1$ to the factors regarded as rational functions  on $X$, say,
$f,h$, can be made regular by multiplying $h$ by an appropriate
polynomial $q\in\C[t]$. Then $\varphi=(f,g):X\to\A^2$, where
$g=qh$, yields a birational morphism. Since every birational
morphism between affine varieties is an affine modification
\cite[Prop.\ 1.1]{KaZa} the claim follows.

For instance the following example from \cite{Be} can be treated
in terms of affine modifications.\erem

\bexa\label{ber} (\cite[Ex. 2.6.1]{Be}, \cite[7.1]{KaZa}) The {\it
Bertin surfaces} are surfaces in $\A^3$ with equations
$$x^ez=x+y^d\,.$$ Every such surface $X$ appears as affine
modification of the plane $\A^2={\rm Spec}\, \C[x,y]$ with center
$(I,\,(x^e))$, where $I=(x^e,x+y^d)\subseteq\C[x,y]$. Actually $X$
is a $\Q$-homology plane with ${\rm Pic}(X)\cong
H_1(X;\Z)\cong\Z/d\Z$, and $\pi=x|X:X\to\A^1$ gives an
$\A^1$-ruling on $X$ with a unique multiple fiber  of multiplicity
$d$ over $x=0$. Whereas the $\A^1_*$-fibration $f=x^{e-1}z : X \to
\A^1$ appears as the orbit map of a $\C^*$-action on $X$ (cf.
Remark \ref{mlinv}.1 below). \eexa

\bprop\label{prop1} Any two
disjoint homology lines $\Gamma_0$ and $\Gamma_1$ on a $\Q$-homology plane $X$
appear as two different fibers of
an $\A^1$-ruling
$\pi:X\to \A^1$. In particular, $X$ arises as in the above construction.
If, moreover,  $X$ is a $\Z$-homology plane then there exists an isomorphism
$X\simeq\A^2$ sending
$\Gamma_0,\Gamma_1$ to two parallel lines.  \eprop

\bproof The second assertion is proven in \cite[\S 9]{Za}.
The first one can be deduced by a similar argument. Namely, ${\rm Pic} (X)$
being a torsion group,
$m\Gamma_0$ is a
principal divisor and so $m\Gamma_0=q^*(0)$ for some $m\in\N$,
$q\in\cO[X]$.
Applying Stein factorisation we
may assume that the general fibers of $q$ are irreducible.
Since $\Gamma_0$ and $\Gamma_1$ are disjoint,
$\Gamma_1$ is
a component of a fiber, say, $F_1=q^{-1}(1)$.
Let the degenerate fibers of $q$
be among the fibers
$F_0=\Gamma_0$, $F_1,\ldots,F_n$, and denote $F$ a general fiber of $q$.
By the Suzuki formula {\it loc. cit.}
$$\sum_{j=0}^n (e(F_j)-e(F))=1-e(F)\,,$$ where all
summands are non-negative by \cite[3.2]{Za}. Since
$e(F_0)=e(\Gamma_0)=1$ we have $e(F_j)=e(F)$ $\forall j=1,\ldots,n$.
It follows by \cite[3.2]{Za} or \cite{Gu}  that either
\begin{enumerate}
[(i)] \item
 the fibers $F_j$ are general
$\forall j=1,\ldots,n$, or
\item
$F\simeq F_j\simeq\A^1_*$ $\forall
j=1,\ldots,n$, or
\item $F\simeq F_j\simeq\A^1$ $\forall
j=1,\ldots,n$.
\end{enumerate}
The case (ii) must be excluded since
$\Gamma_1\subseteq F_1$ and
$e(\Gamma_1)=1$.
If (i) holds then $F_1\simeq\pi^*(1)$ is a general fiber of $q$,
hence $F\simeq F_1\simeq\A^1$.
Thus in any case
all fibers of $\pi$ are isomorphic to $\A^1$ and so,
 $\Gamma_0,\Gamma_1$ are fibers of the $\A^1$-ruling
$\pi=q:X\to\A^1$.
\eproof

Now one can easily deduce Theorem \ref{theo}(a).

\smallskip

{\em Proof of Theorem \ref{theo}(a)}. $X\setminus \Gamma$ being a
smooth affine surface with $\bar k (X\setminus \Gamma)=-\infty$,
there exists an $\A^1$-ruling on $X\setminus \Gamma$
\cite[2.1.1]{Mi}. The curve $\G_0=\G$ and a general fiber, say,
$\Gamma_1$ of this ruling provide two disjoint homology lines on
$X$. By Proposition \ref{prop1} $\Gamma_0$ and $\Gamma_1$ are  two
different fibers of an $\A^1$-ruling $\pi:X\to\A^1$ on $X$ and so,
$\Gamma=\Gamma_0$ is an affine line stable under an effective
$\C_+$-action on $X$ along this ruling, see e.g.
\cite[1.6]{FlZa3}. Now the assertion follows easily. \qed

\section{$\Q$-homology planes with a $\C^*$-action}

To deduce Theorem \ref{theo}(b) we recall first Example 1 in
\cite{KiKo}, cf. also Examples 3.8, 3.9 in \cite{GMMR}.

\smallskip

\noindent {\bf Construction.} The construction begins with a
divisor $D_0=M_a+\bar M_a+F_0+F_1+F_\infty$ on a Hirzebruch
surface $\Sigma_a$, where $F_0,F_1,F_\infty$ are 3 distinct fibers
of the standard projection $\pi_0:\Sigma_a\to\PP^1$ and
$M_a,\,\bar M_a$ are two disjoint sections with $M_a^2=-a$. We may
suppose that $F_j=\pi_0^{-1}(j)$ ($j=0,1,\infty\in\PP^1$).
Besides, the construction  involves a sequence of {\it inner}
blowups $\mu:V\to \Sigma_a$ over $D_0$ i.e., successive blowups
with centers at double points on $D_0$ or on its total transforms.
This results in a $\Q$-homology plane $X=V\setminus D$, where
$D\subseteq \mu^{-1}(D_0)$ is a suitable SNC tree of rational
curves on $V$ (see a description below). The induced
$\PP^1$-ruling $\bar\pi : V \to \PP^1$ restricts to an untwisted
$\A^1_*$-fibration $\pi:X\to \A^1$.

More precisely $\mu$ replaces the fiber $F_j$ ($j=0,1$)
by a linear chain of smooth rational curves with a
unique $(-1)$-curve $E_j$,
which is a multiple component of the corresponding
divisor $\bar\pi^*(j)$.
The dual graph of the chain $\bar\pi^{-1}(0)$ has a
sequence of weights
$[[-n,-1, \underbrace{-2,\ldots,-2}_{n-1}]]$, where for
the strict transform $F'_0$ of $F_0$ on $V$ one has
${F'}_0^2=-n\le -2$.
The boundary divisor $D$ appears as the total transform
of $D_0$ in $V$ with the components $E_0,E_1$
and $F_0'$ being deleted. In the affine part $X=V\setminus D$
these deleted components
form the only degenerate fibers of $\pi$,
namely $\pi^*(0)=\G+n(E_0\cap X)$ and $\pi^*(1)=m(E_1\cap X)$,
where $\G= F_0'\cap X$
and $m,n\ge 2$.
Thus the only reducible affine fiber
$\pi^{-1}(0)=(F_0'+E_0)\cap X$ is isomorphic to the cross
$\A^1\wedge\A^1=\{xy=0\}\subset\A^2$. Furthermore
$\pi^{-1}(1)= E_1\cap X\simeq \A^1_*$ is
an irreducible multiple fiber. Finally
$\bar\pi^{-1}(F_\infty)=F_\infty'\subseteq D$. A computation in
\cite[Example 1]{KiKo} shows that $\bar k (X\setminus \G)=0$ if $m=n=2$
and
$\bar k (X\setminus \G)=1$
otherwise.

\brem\label{cuch} One could consult e.g. \cite{FlZa1} for a construction
giving all $\Q$-homology planes $X$
with an $\A^1_*$-fibration $\pi:X\to B$. In the terminology of
\cite{FlZa1}, such a surface
with a twisted (untwisted) $\A^1_*$-fibration over $B=\A^1$,
 $B=\PP^1$, respectively,
 is said to be of type $A1$,  $A2$ ($B1$, $B2$), respectively.
Thus the surface $X$
as in the Construction above is
of type $B1$, with
a comb attachment applied at $F_0$ and with $F_1$ replaced by a fiber of
a {\it broken chain} type in the terminology of \cite{FlZa1}. \erem

In the following lemma we prove the first assertion
of Theorem \ref{theo}(b).
We recall that a $\C^*$-action on $X$ is {\it hyperbolic}
if its general orbits are
closed, {\it elliptic} if it
possesses an attractive or repelling fixed point in $X$,
and {\it parabolic} if its fixed point set is one-dimensional.

\blem\label{conac} If $(X,\G)$ satisfies
the assumptions of Theorem
\ref{theo}(b)
then $\G$ is an orbit closure of an effective hyperbolic
$\C^*$-action on $X$. \elem

\bproof
According to Theorem 1.1 in \cite{KiKo}
(cf. Theorem 3.10 in \cite{GMMR}),
under our assumptions
$(X, \Gamma)$ is one of the pairs as
in the  above Construction.
There exists an effective
$\C^*$-action
on $\Sigma_a$ along the fibers of
$\pi_0$ with the
fixed point set equal to $M_a\cup\bar{M_a}$.
By induction on the number of blowups this
$\C^*$-action lifts to $V$ stabilizing
the total transform $\mu^{-1}(D_0)$.
Indeed the centers of successive inner
blowups in $\mu$ are fixed under the $\C^*$-action
constructed on the
previous step, and so Lemma 2.2(b) in
\cite{FKZ} applies. It follows that the curve
$D\subseteq \mu^{-1}(D_0)$ as
in the Construction above
is stable under the lifted $\C^*$-action,
so this action restricts to a hyperbolic
$\C^*$-action on $X=V\setminus D$.
In turn, the affine line
$\G$ on $X$ as in the Construction
is an orbit
closure for this restricted action, as required.
\eproof

The resulting surface $X$ with a hyperbolic
$\C^*$-action admits the following
description in terms of
the DPD orbifold presentation\footnote{i.e.
the Dolgachev-Pinkham-Demazure presentation.} as elaborated
in \cite{FlZa2}.

\smallskip

\noindent {\bf DPD presentation}. Let $C={\rm Spec}\, A_0$ be a
smooth affine curve and $(D_+,D_-)$ be a pair of $\Q$-divisors on
$C$ with $D_++D_-\le 0$. Letting $A_{\pm k}=H^0(C, \cO(kD_\pm))$,
$k\ge 0$ we consider the graded $A_0$-algebra
$A=A_0[D_+,D_-]=\bigoplus_{k\in\Z} A_k$ and the associated normal
affine surface $X={\rm Spec}\, A$. The grading determines,  in a
usual way, a graded semisimple Euler derivation $\delta$ on $A$,
where $\delta (a_k)=ka_k$ $\forall a_k\in A_k$, and, in turn, an
effective hyperbolic $\C^*$-action on $X$. Vice versa, any
effective hyperbolic $\C^*$-action on a normal affine surface with
the orbit space $C$ arises in this way \cite[4.3]{FlZa2}.

Let $\pi:X\to C$ be the orbit map. Given
a point $p\in C$ we let $m_\pm(p)$
denote the minimal positive integer
such that $m_\pm(p) D_\pm(p)\in\Z$.
In case where  $(D_++D_-)(p)=0$ we set $m(p)=m_\pm(p)$.
If $(D_++D_-)(p)<0$ then
the fiber $\pi^{-1}(p)$ is reducible, isomorphic to the cross
$\A^1\wedge\A^1$ in $\A^2$ and consists of two
orbit closures $\bar O_p^\pm$. Its unique double point $p'$
is a fixed point; $p'$ is smooth on $X$
if and only if $(D_++D_-)(p)=-1/m_+m_-$ \cite[4.15]{FlZa2}.
Actually $m_\pm(p)$ are the multiplicities
of the curves $\bar O_p^\pm$, respectively, in the divisor $\pi^*(p)$.

In case where $(D_++D_-)(p)=0$
the fiber $O_p=\pi^{-1}(p)\simeq\A^1_*$
is irreducible of multiplicity
$m(p)$ in $\pi^*(p)$.

The inversion $\lambda\longmapsto\lambda^{-1}$ in $\C^*$ results
in interchanging $D_+$ and $D_-$, respectively, $\bar O_p^+$ and
$\bar O_p^-$. Passing from the pair $(D_+,D_-)$ to another one
$(D'_+,D'_-)=(D_++D_0,D_--D_0)$ with a principal divisor $D_0$ on
$C$ results in passing from $A$ to an isomorphic graded
$A_0$-algebra $A'$, so the corresponding $\C^*$-surfaces  are
equivariantly isomorphic over $C$.

\blem\label{qacyc} Given
a normal affine surface $X={\rm Spec}\, A$ with a hyperbolic $\C^*$-action
determined by a pair $(D_+,D_-)$ of $\Q$-divisors
on the affine
curve $C={\rm Spec}\, A_0$  with $D_++D_-\le 0$, we
denote by $p_1,\ldots,p_l$, $q_1,\ldots,q_k$
the points of $C$ with
$(D_++D_-)(p_j)<0$,
$(D_++D_-)(q_i)=0$ and $m(q_j)\ge 2$, respectively.
Letting $\pi:X\to C$ be
the orbit map we assume that $C\simeq\A^1$ and that $X$ is smooth that is,
$(D_++D_-)(p_j)=-1/m_+(p_j)m_-(p_j)$
$\forall j=1,\ldots,l$. Then the following hold.
\begin{enumerate} [(a)] \item $e(X)=l$.
\item
${\rm Pic} (X)\otimes\Q=0$ if and only if $l\le 1$ that is, $\pi$
has at most one reducible fiber.
\item
Moreover $X$ is $\Q$-acyclic if and only if $l=1$. In the latter case
$\pi:X\to \A^1$ is an untwisted $\A^1_*$-fibration\footnote{It is of type B1
in the classification of \cite{FlZa1}.}.
\end{enumerate}\elem

\bproof (a) holds by the additivity of the Euler characteristic, and (b)
follows from the description of the Picard group
${\rm Pic} (X)$ in \cite[4.24]{FlZa2}.
For a smooth rational affine surface $X$
we have $b_3=b_4=0$ and $b_1=\rho (X)$,
where $\rho (X)$
is the Picard number
of $X$ \cite[Ch. 3, 4.2.1]{Mi}. Thus $X$ is $\Q$-acyclic
if and only if $e(X)=1$ and
${\rm Pic} (X)\otimes\Q=0$,
whence
(c) follows.
\eproof

\blem\label{dpd} Every $\Q$-homology plane $X$ as in the above
Construction is isomorphic to a $\C^*$-surface ${\rm Spec}\,
A_0[D_+,D_-]$, where $A_0=\C[t]$ and $$D_+=\frac{e}{m}[1],\quad
D_-= -\frac{1}{n}[0]-\frac{e}{m}[1]\,\,\,\mbox{with}\,\,\,
0<e<m,\,\gcd (e,m)=1,\,\,m,\,n\ge 2\,.$$ Conversely, every
$\C^*$-surface with such a DPD-presentation appears via the above
Construction. \elem

\bproof
By our Construction, the degenerate fibers of the
induced $\A^1_*$-family $\pi :X \to \A^1$ are
$\pi^*(0)=n(E_0\cap X)+\G$ and $\pi^*(1)=m(E_1\cap X)$,
where $\G=F_0'\cap X$, $E_0\cdot \G=1$
and $E_j$ (j=0,1)
is the unique $(-1)$-curve in the fiber $\bar\pi^{-1}(j)$
of the induced $\PP^1$-ruling
$\bar\pi:V\to \PP^1$. Clearly, all these curves are orbit closures for
the $\C^*$-action on $X$ as in Lemma \ref{conac}.
We may suppose that, in the notation as above,
$\G=O_0^+,\,\, E_0\cap X=O_0^-$ and $E_1\cap X=O_1$ so that
$$k=l=1,\,\,\, p_1=0,\,\,\, q_1=1,
\,\,\, m_+(0)=1,\,\,\, m_-(0)=n\,\,\,\mbox{and}\,\,\, m(1)=m\,.$$
Since every integral divisor on $C=\A^1$ is principal, passing to
an equivalent pair of $\Q$-divisors we may achieve that
$(D_+,D_-)$ is a pair as in the lemma. This proves the first
assertion. The converse easily follows by virtue of Lemma
\ref{qacyc}. \eproof

\brems\label{mlinv} 1. According to \cite[5.5]{FlZa3}, for $e=1$
and $m\mid n$ the above surfaces actually coincide with the Bertin
surfaces from Example \ref{ber}.

2. The formula for the canonical divisor
in \cite[4.25]{FlZa2} gives in our case
$$K_X=-(e(n-1)+1)[O_1],\qquad\mbox{where}\quad m[O_1]=0\,.$$
Therefore $K_X=0$ if and only if $e(n-1)\equiv -1\, {\rm mod}\,
m$. The question arises whether, among the $\C^*$-surfaces from
Lemma \ref{dpd} satisfying the latter condition, the Bertin
surfaces are the only hypersurfaces.

3. For $n>1$ the fractional part $\{D_-\}$ in Lemma \ref{dpd}
is supported on 2
points, hence by \cite[4.5]{FlZa3}
the surface $X$ as in Lemma \ref{dpd} admits a unique
$\A^1$-ruling $X\to\A^1$ (i.e., $X$ is of class ML$_1$
in the terminology of \cite{GMMR}).

In contrast, for
$n=1$ there is a second $\A^1$-ruling $X\to\A^1$,
so $X$ has trivial
Makar-Limanov invariant. In particular for $e/m=1/2$ and $n=1$
by virtue of \cite[5.1]{FlZa3},  $X\simeq
\PP^2\setminus\Delta$, where $\Delta$ is a smooth conic in $\PP^2$.

4. Following \cite[4.8]{FlZa2} it is possible to define, by
explicit equations, a family of surfaces in $\A^4$, not
necessarily complete intersections, whose normalizations are the
$\Q$-homology planes in the above Construction. \erems

Now we are ready to complete the proof of Theorem \ref{theo}(b).

\smallskip

{\em Proof of Theorem \ref{theo}(b)}.
Since the fractional part $\{D_+\}$ of the divisor $D_+$
as in Lemma \ref{dpd}
is supported on one point,
 there exists a
graded locally nilpotent derivation on $A$ of positive degree,
see \cite[2.2, 3.23]{FlZa3}.
This  derivation
generates an effective $\C_+$-action on $X$, and also an action
of a semidirect product
$G=\C^* \ltimes\C_+$ with
an open orbit $U\simeq\A^1\times\A^1_*$.
Moreover by \cite[3.25]{FlZa3}, the orbit map $X\to\A^1$
of the associate $\C_+$-action has a unique irreducible multiple
fiber $\G'=\bar O^-_0$ ($=E_0\cap X$) of multiplicity $m\ge 2$.
General orbits of this $\C_+$-action on $X$ being transversal
to $\G$, the action moves $\G$, as stated. \qed

\section{Isotrivial families of curves and $\C^*$-actions}
To indicate a proof
 of Theorem \ref{theo}(c) let us
recall first a necessary result from \cite{LiZa, Za}. For the sake
of completeness we sketch the proof.

\blem\label{iso} {\rm (\cite[Lemma 5]{LiZa})} Let $X^*$ be a
smooth affine surface and $\pi:X^*\to \A^1_*$ be a family of
curves without degenerate fibers which is not a twisted
$\A^1_*$-family. Then $\pi$ is equivariant with respect to a
suitable effective $\C^*$-action on $X^*$ and a nontrivial
$\C^*$-action on $\A^1_*$. \elem

\bproof Let $F$ denote a general fiber of $\pi$.
In the case where $F\simeq \A^1$
the surface $X^*$ admits a completion
which is a Hirzebruch surface $\Sigma_a$ with the boundary divisor
$D=\Sigma_a\setminus X^*$
consisting of a section and two fibers. It follows that $\pi$ is a
trivial family, which implies the assertion.
The same argument applies if $F\simeq\A^1_*$
since in this case by our assumption $\pi$ is untwisted.

Suppose further that $e(F)<0$ i.e. that $F$ is a hyperbolic curve.
By Bers' Theorem the Teichmuller space corresponding to $F$, with
its natural complex structure, is biholomorphic to a bounded
domain in $\C^M$ for some $M>0$, hence is as well hyperbolic.
Therefore the family $\pi$ over a non-hyperbolic base $\A^1_*$ is
isotrivial i.e., its fibers are all pairwise isomorphic. Since
${\rm Aut} (F)$ is a finite group the monodromy $\mu\in {\rm Aut}
(F)$ of the family $\pi$ has finite order, say, $N$. After a
cyclic \'etale base change $z\longmapsto z^N$ we obtain a trivial
family $F\times \A^1_*\to\A^1_*$, which is a cyclic \'etale
covering of the given family $\pi$. The standard $\C^*$-action on
its base lifts to a free $\C^*$-action on $F\times \A^1_*$
commuting with the monodromy $\Z/N\Z$-action. Therefore the lifted
$\C^*$-action descends to $X^*$ so that $\pi$ becomes equivariant
with respect to the $\C^*$-action $\lambda . z=\lambda^N z$ on
$\A^1_*$, as needed. \eproof

\brem\label{twi} The $\A^1_*$-family of orbits of a
hyperbolic $\C^*$-action on an affine surface is always untwisted
\cite{FKZ}.
Hence the conclusion of Lemma \ref{iso} does not hold
for twisted $\A^1_*$-families.
\erem

{\em Proof of Theorem \ref{theo}(c)}. Let $\G$ be a non-smooth
homology line on a $\Q$-homology plane $X$, and let $m\G=f^*(0)$
for a suitable $m\in\N$ and a primitive regular function
$f\in\cO(X)$ with irreducible general fiber $F$ (cf. the proof of
Proposition \ref{prop1}). Let $p'\in\G$ be a singular point of
$\G$ with Milnor number $\mu>0$. In a suitable small spherical
neighborhood $B$ of $p'$, the function $f^{1/m}$ is holomorphic
and its general fiber say $R$ (which is the Milnor fiber of
$(\G,p')$) is a Riemann surface with boundary of positive genus
$g=\mu/2$ \cite[10.2]{Mil}. For a fixed general fiber $F$ of $f$
sufficiently close to $\G$, the intersection $F\cap B$ is a
disjoint union of $m$ copies of the Milnor fiber $R$, hence $F$ as
well is of positive genus.

Therefore $e(F)<0$. By Lemma 3.2 in \cite{Za}, since
$e(X\setminus\G)=0$ the family $\pi=f\vert (X\setminus\G) :
X\setminus\G\to \A^1_*$ has no degenerate fiber, and so Lemma
\ref{iso} applies.

As a matter of fact,
the $\C^*$-action on $X\setminus\G$
as in Lemma \ref{iso} extends to an elliptic
$\C^*$-action on
$X$ making $f$ equivariant and $p'$ an
attractive or repelling fixed point.
For a $\Z$-homology plane $X$,
the existence of such an extension
was shown
in \cite{LiZa} and in \cite{Za} in two
different ways. The both proofs
work {\it mutatis mutandis} in our more general setting.
We choose below to follow the lines of
the proof of Lemma 6 in \cite{LiZa}.

Let $\bar F$ be a smooth projective model of $F$.
The cyclic \'etale covering
$\rho : F\times \A^1_*\to X\setminus\G$
as in the
proof of Lemma \ref{iso} extends to an equivariant rational map
which fits into the commutative
diagram \be
\label{diagramm}
\begin{diagram}
\bar F\times \PP^1 & \rDashto^{\rho} & V \\
\dTo_{\rm pr_2} & &\dTo_{\bar\pi}\\
\PP^1& \rTo^{z\longmapsto z^N} & \PP^1
\end{diagram}
\ee where $V$ is a smooth equivariant SNC completion of
$X\setminus\G$. If the $\C^*$-action on $X\setminus\G$ possesses
an orbit $O$ which is not closed in $X$ then so are all orbits and
the action extends to $X$. Indeed the closure $\bar O$ meets $\G$
in one point say $q$, and $n\bar O=h^*(0)$ for some regular
function $h$ on $X$ and for some $n\in\N$. Let $P$ be the
connected component of the polyhedron $|f|\le
\varepsilon,\,\,|h|\le \varepsilon$ which contains $q$. Since
$\lambda . h=\lambda^k \cdot h$ for some $k\in\Z$, for every
$\lambda\in\C^*$ with $|\lambda|=1$ the complement $P\setminus\G$
is stable under the action of $\lambda$. The polyhedron $P$ being
compact for a sufficiently small $\varepsilon>0$, such an action
on $P\setminus\G$ extends across $\G$. Hence it also extends
through $\G$ to an action on the whole $X$ for all
$\lambda\in\C^*$ with $|\lambda|=1$, and then also for all
$\lambda\in\C^*$.

The indeterminacy set of $\rho$ being at most $0$-dimensional,
$\rho$ restricts
to the fiber over $0\in\PP^1$ yielding a morphism
$\rho : \bar F\to\bar\pi^{-1}(0)$. The following alternative holds:
either $\rho (\bar F)=p\in\G$, or $\rho (\bar F)=\bar\G$,
or finally
$\rho (\bar F)\cap \G=\emptyset$. Let us show that
the last two possibilities cannot occur.

Indeed supposing that $\rho (\bar F)=p\in\G$
($\rho (\bar F)=\bar\G$,
respectively)
the general orbits
of the $\C^*$-action on $X\setminus \G$ are not closed in $X$,
and the action
extends to an elliptic (parabolic, respectively)
$\C^*$-action on $X$.
Since the fixed
point set of a parabolic $\C^*$-action on a normal
affine surface is smooth (see \cite[\S 3]{FlZa2}),
the
latter case must be excluded.

To exclude the last possibility, suppose on the contrary that
$\rho (\bar F)\cap \G=\emptyset$.
Letting $p_1,\ldots,p_k\in\bar F\times \{0\}$
be the indeterminacy points of
$\rho$ on the central fiber, we observe that under our assumption,
all $\C^*$-orbits $\rho
(\{p\}\times\A^1_*)$ are closed in $X$, because general orbits are.
The orbits
$\rho (\{p_i\}\times\A^1_*)$, $i=1,\ldots,k$, meet any fiber
$F_\xi=f^{-1}(\xi)$, $\xi\in\A^1$, in a finite set, say, $T$.

Fixing further a general fiber $F=F_\xi$ sufficiently close to
$\G$ and a
sufficiently small neighborhood $\omega$ of the finite set
$T\cup (\bar F\setminus F)$
in $\bar F$, we let $K=\bar F\setminus \omega\subseteq F$.
Under our assumptions $K$ is a compact Riemann surface of
positive genus with boundary, and $B\cap \lambda . K=\emptyset$ for all
sufficiently small $\lambda\in\C^*$. Hence $F\cap \lambda^{-1} . B
\subseteq\omega$ is a disjoint union of Riemann surfaces of genus $0$.
On the other hand
$$F\cap \lambda^{-1} . B \cong B\cap \lambda . F = B\cap
F_{\lambda^N\xi}$$ is a disjoint union of
$m$ copies
of the Milnor fiber $R$ of the analytic
plane curve singularity $(\G,p')$.
This is a
contradiction because $R$ is of positive genus.

Thus $\G$ is stable under the extended elliptic $\C^*$-action on
$X$. So  $\G$ is an orbit closure of this action and the singular
point $p'\in\G$ is a fixed point of the action. Consider an
equivariant embedding $X\hookrightarrow \A^N$ which sends $p'$ to
the origin, where $\A^N$ is equipped with a linear $\C^*$-action,
and fix an equivariant linear projection $\A^N\to T$, where
$T\simeq\A^2$ is the tangent plane of $X$ at $p'=\bar 0$. This
projection  restricted to $X$ gives an equivariant isomorphism
$X\simeq\A^2$, where the $\C^*$-action on $\A^2$ is linear
(indeed, both actions have the origin as an attractive fixed
point). In appropriate linear coordinates the latter linear action
is diagonal:
 $\lambda .
(x,y)\longmapsto (\lambda^lx,\lambda^ky)$ with $\gcd (k,l)=1$. So
either the image of $\G$ is one of the axes,
which contradicts the assumption that $\G$ is
singular,
or it is a curve
$\alpha x^k-\beta y^l=0$ for some $\alpha,\beta\in \C^*$.
This proves (c) of Theorem \ref{theo}. Now the proof of Theorem \ref{theo} is
completed.
\qed

\end{document}